\documentclass{article}

\usepackage[margin=1in]{geometry}
\usepackage[all]{xy}
\usepackage{amsmath,amssymb,amsthm}
\usepackage{chngcntr}
\usepackage{color}
\usepackage{comment}
\usepackage{enumerate}
\usepackage{graphicx}
\usepackage{mathrsfs}
\usepackage{mathtools}
\usepackage{paralist}
\usepackage{amscd}
\usepackage[OT2,T1]{fontenc}
\DeclareSymbolFont{cyrletters}{OT2}{wncyr}{m}{n}
\DeclareMathSymbol{\Sha}{\mathalpha}{cyrletters}{"58}

\theoremstyle{definition}

\newtheorem{theorem}{Theorem}
\newtheorem*{theorem*}{Theorem}
\newtheorem{proposition}[theorem]{Proposition}
\newtheorem*{proposition*}{Proposition}
\newtheorem{lemma}[theorem]{Lemma}
\newtheorem*{lemma*}{Lemma}

\newtheorem*{corollary*}{Corollary}

\newtheorem*{example*}{Example}

\newtheorem*{definition*}{Definition}
\newtheorem{remark}[theorem]{Remark}
\newtheorem*{remark*}{Remark}

\newtheorem*{fact*}{Fact}

\newtheorem*{observation*}{Obzervation}

\newtheorem*{claim*}{Claim}
\newtheorem{conjecture}[theorem]{Conjecture}

\newtheorem*{notation*}{Notation}

\newtheorem*{convention*}{Convention}

\newtheorem*{assumption*}{Assumption}

\numberwithin{theorem}{section} 	 
\numberwithin{equation}{section}	 

\mathtoolsset{showonlyrefs=true} 

\newcommand{\Z}{\mathbb{Z}}
\newcommand{\Q}{\mathbb{Q}}

\newcommand{\fr}[2]{\frac{#1}{#2}}

\newcommand{\ol}[1]{\overline{#1}}

\begin{document}

\title{ON SOME VALUES WHICH DO NOT BELONG TO THE IMAGE OF RAMANUJAN'S TAU-FUNCTION}
\author{AKIHIRO GOTO}
\date{}
\maketitle

\begin{abstract}
Lehmer conjectured that Ramanujan's tau function never vanishes.
As a variation of this conjecture, it is proved that 
\begin{equation*}
\tau(n)\neq \pm \ell, \pm 2\ell, \pm 2\ell^2,
\end{equation*}
where $\ell<100$ is an odd prime, by Balakrishnan, Ono, Craig, Tsai and many people.
We have proved that
\begin{equation*}
\tau(n)\neq \pm \ell, \pm 2\ell, \pm 4\ell, \pm 8\ell
\end{equation*}
for any $n\geq 1$ except 14 cases, where $\ell<1000$ is an odd prime.
\end{abstract}

\tableofcontents


\addcontentsline{toc}{section}{Introduction and main results}
\section*{Introduction and main results}
For ${\rm Im }(z)>0$, we define Ramanujan's tau function
\begin{equation*}
\Delta(z):=\sum_{n=1}^{\infty}\tau(n)q^n:=q\prod_{n=1}^{\infty}(1-q^n)^{24}, q:=e^{2\pi iz}.
\end{equation*}
Then $\Delta$ is the normalized Hecke eigenform with weight 12 and level 1.

In 1947, Lehmer conjectured that $\tau$ is never vanish \cite{L}.
This conjecture is still unproven.
Serre showed that the natural density of the set of primes $p$ which $\tau(p)=0$ within primes is 0.

Naturaly, as a variation of Lehmer's conjecture, it can be considered that the existence of nonzero integer which never appears as $\tau$-values.
Murty, Murty and Shorey showed in \cite{MMS} that 
for any odd integer $a$, $\tau(n)=a$ holds for at most finitely many $n$, i.e., the Lang--Trotter conjecture {\cite{LT} is true if $a$ is odd.
In 2013, Lygeros and Rozier \cite{LR}, showed $\tau(n)\neq \pm 1$ for any $n>1$.
In 2020, Balakrishnan, Ono, Craig, Tsai obtained that
$$
\pm 3, \pm 5, \pm 7, \pm 13, \pm 17, -19, \pm 23, \pm 37, \pm 691
$$
are not $\tau$-values.
Hanada and Madhukara \cite{HM} proved that 
$$
-9, \pm 15, \pm 21, -25, -27, -33, \pm 35, \pm 45, \pm 49, -55, \pm 63, \pm 77, -81, \pm 91 
$$
are not $\tau$-values also.
Dembner and Jain \cite{DJ} showed that
\begin{equation} \label{tauneqpm1prime691}
\tau(n)\not \in \{\pm \ell : \ell <100, \text{ odd prime} \} \cup \{\pm 5^m : m\geq 0\}
\end{equation}
for any $n>1$.
Bennett, Gherga, Patel and Siksek proved that for $n>1$,
\begin{equation} 
\tau(n)\neq \pm \ell^a, 3^b5^c7^d11^e, \text{ where $\ell$ is odd prime $<100$}, a, b, c, d, e\geq 0 
\end{equation}
holds in \cite[Corollary 1.1 and Theorem 6]{BGPS}.
Lin and Ma showed that if $\tau(n)$ is odd for $n>3$, 
\begin{equation} \label{lemfromLMcor1.3}
|\tau(n)|>|\tau(3)|=252 
\end{equation}
holds in \cite[Corollary 1.3]{LM}.
Above results on odd $\tau$-values were showed by using the following fact:
\begin{equation} \label{tauodd}
\Delta(z)\equiv \sum_{n=0}^{\infty}q^{(2n+1)^2} \bmod 2,
\end{equation}
which can be proved applying the Jacobi's triple product identity.

On the other hand, as even cases, Balakrishnan, Ono, Tsai showed that if $3 \leq \ell <100$ is odd prime, then
\begin{equation} \label{EVTcor1.2}
\tau(n)\not \in \{\pm 2\ell \}\cup \{\pm 2\ell^2\}\cup\{\pm 2\ell^3:\ell\neq 59\}
\end{equation}
for any $n\geq 1$ in \cite{EVT}.
As general theory of the number of prime divisors of $\tau$-values, there is the following theorem:
\begin{theorem}[\cite{BCKT}]
We define $\omega(n) (\text{resp. }\Omega(n))$ as the number of prime divisors of $n$ (resp. with multiplicity).
For $n \geq 1$,
$$
\Omega(\tau(n))\geq \sum_{i=1}^r(\sigma_0(d_i)-1) \geq \omega(n)
$$
holds, where $n=p_1^{d_1-1}\cdots p_r^{d_r-1}$ and each $p_i$ is ordinary prime, which means that $p_i\nmid \tau(p_i)$.
\end{theorem}


\begin{remark}
The inequation $\Omega(\tau(n))\geq \omega(n)$ holds even if $n\geq 1$ is divisible by some non-ordinary prime $p$, since for any prime power $q$, we have $|\tau(q)|>1$(by \eqref{tauneqpm1prime691}) and the multiplicativity of $\tau$ (Theorem \ref{RCj} \eqref{RC1}).
\end{remark}


In this paper, we were inspired by the above results, we study more numbers which do not belong to the image of the $\tau$-function.
To state the main theorems, we set
\begin{gather}
L_{1}^+:=\{461\}, L_{1}^-:=\{599\}, L_2^{+}:=\emptyset, L_2^{-}:=\{587\}, \\
L_{4}^+:=\{23, 449, 569, 863\}, L_{4}^-:=\{241, 397, 811\}, L_{8}^+:=\{457\}, L_{8}^-:=\{3, 293, 983\}.
\end{gather}
Let $\ell<1000$ be an odd prime and $\varepsilon=\pm 1$.
Then, the following two theorems are main results:
\begin{theorem} \label{mainthmpmell}
If $\tau(n)=\varepsilon \ell$, we must have $n=p^4$, where $p$ is prime and $\ell\in L_1^{\varepsilon}.$
\end{theorem}


\begin{theorem} \label{mainthm248}
Let $t=2, 4, 8$.
If 
$$
\tau(n)=\varepsilon t\ell,
$$
then the one of the following is true:
\begin{enumerate}[(1)]
\item \label{fcmt}
There exists prime $p \mid n$ such that $\tau(p)=2$ and $\gcd(p, n/p)=1$,

\item \label{scmt}
$n$ is prime and $\ell \in L_{t}^{\varepsilon}$.
\end{enumerate}
\end{theorem}


For example, if $\tau(n)=-4\ell$, by using distinct primes $p_1, p_2, p_3$, we can write
$$
n=
\begin{cases}
p_1p_2p_3^4 & \text{such that } \tau(p_1)=\tau(p_2)=2, \tau(p_3^4)=-599, \\
p_1p_2 & \text{such that } \tau(p_1)=2, \tau(p_2)=-2 \cdot 587, \\
p_1 & \text{such that } \ell=241, 397, 811.
\end{cases}
$$


\begin{remark}
\begin{enumerate}[(1)]
\item
In each case, some congruent condition of $p$ must be satisfied.
For example, if $\tau(p)=-2\cdot 587$, then $p$ must satisfy 
$$
p\equiv 
\begin{cases}
1 & \bmod 3, \\
1 & \bmod 4, \\
2 & \bmod 5, \\
1 & \bmod 7.
\end{cases}
$$
Moreover, $(p/23)=1$ and $p$ cannot be written as $a^2+23b^2, a,b \in \Z$.
See also \eqref{RC} or more convenient Lemma \ref{devencongruences}.
\end{enumerate}
\end{remark}


Among the values $\varepsilon t \ell$ in Theorem \ref{mainthm248}, only $-8\cdot 3(=\tau(2))$ was confirmed to actually appear.


The first example of $\tau(n)=\pm 2\ell, \pm 4 \ell, \pm8\ell, \ell:\text{odd prime}$:
$$
\tau(277) =-2\cdot 8209466002937, \tau(1297) = 2\cdot 58734858143062873,
$$
$$
\tau(163)=-4\cdot 89458189897, \tau(4603)=4\cdot 56958468932026008713,
$$
$$
\tau(2)=-8\cdot 3, \tau(967)=8\cdot 2311913038549627, \tau(2647)=8\cdot 1344910678663379137.
$$


The proof of main results is a combination of tools which is used in previous researches \cite{BCO, BCKT, EVT, DJ}.
Roughly speaking, thanks to Proposition \ref{BCKT}, we may consider only $\tau$-value of prime powers.

In \cite{Lucas}, for any Lucas sequence $\{u_n\}$, it is completely known about $u_n$ which never has primitive prime divisors (the definition is the subsection \ref{Ls}).
In particular, they showed that for all Lucas sequence $\{u_n\}_{n=1}^{\infty}$, if $n>30$, $u_n$ has a primitive prime divisor.
Therefore, if $n>30$,
\begin{equation}
u_n\neq \prod_{\text{prime }\ell \mid u_1\cdots u_{n-1}}\ell^{a_{\ell}}, a_{\ell}\geq 0
\end{equation}
holds.
Applying these facts, we proved Lemma \ref{methodEVT}:
\begin{lemma*}[Lemma \ref{methodEVT}]
Let $p, \ell$ be primes, $\ell \neq 2, d\geq 1$.
If
$$
\tau(p^{d-1})=\pm 2\ell, \pm 4\ell, \pm 8\ell,
$$
we have 
$$
d=
\begin{cases}
2, & \tau(p^{d-1})=\pm 2 \ell, \pm 8\ell, \\
2, 4, & \tau(p^{d-1})=\pm 4\ell.
\end{cases}
$$
Moreover, we have that $\{\tau(p^{i-1})\}_{i=1}^{\infty}$ is a Lucas sequence.
If $\tau(p^3)=\pm 4\ell$, then we have $\tau(p)=2$ and that $\ell$ is primitive prime divisor of $\tau(p^3)$.
\end{lemma*}

By using Ramanujan's congruence \eqref{RC} also, eventually, we may consider only $\tau(p), \tau(p^3), \tau(p^4)$.
Finally, we apply Thue equation or Dembner and Jain's method \cite{DJ}, for each individual case.

We will remark on the relation between main results and the Atkin--Serre conjecture:
\begin{conjecture}[Atkin--Serre \cite{Serre}] 
Let $f\in S^{\text{new}}_k(\Gamma_0(q))$ be a non-CM newform of weight $k\geq 4$.
For each $\varepsilon$, there exist constants $c_{\varepsilon, f}>0$ and $c'_{\varepsilon, f}>0$ such that if $p>c'_{\varepsilon, f}$, then 
\begin{equation} \label{asconj}
|a_f(p)|\geq c_{\varepsilon, f}p^{\fr{k-3}{2}-\varepsilon}.
\end{equation}
\end{conjecture}

Deligne showed that if $f$ is a normalized Hecke eigenform of even weight $k\geq 2$, there exists $0\leq \theta_p\leq \pi$ such that $a_f(p)=2p^{(k-1)/2}\cos \theta_p$.

Sato--Tate conjecture which was proved by Barnett-Lamb, Gehrarty, Harris, and Taylor \cite{STS1, STS2} states that for any $0\leq \alpha<\beta\leq \pi$, we have 
\begin{equation} \label{ST}
\displaystyle \lim_{x\to \infty} \frac{\#\{p \leq x : \alpha \leq \theta_p\leq \beta\}}{\#\{p \leq x\}}=\frac{2}{\pi}\displaystyle \int_{\alpha}^{\beta}\sin^2\theta d\theta.
\end{equation}
Recently, Newton--Thorne \cite{NT1} proved that the $n$-th symmetric power $L$-function of $f$ is the $L$-function of an automorphic representation of ${\rm GL}_n(\mathbb{A}_{\Q})$ for any $n\geq 1$, where $\mathbb{A}_{\Q}$ is the ring of Ad\`{e}les over $\Q$.
By using this breakthrough, Thorner proved a strong version of \eqref{ST} with an effective error term.

In \cite{GTW}, Gafni, Thorner, Wong proved that \eqref{asconj} is true except for a density zero set of primes using an effective Sato--Tate theorem \cite{Thorner}.
More precisely, they showed that except for a density zero set $S(f)$ of primes,
\begin{equation}
2p^{\fr{k-1}{2}}\cdot \fr{\log \log p}{\sqrt{\log p}}<|a_f(p)|.
\end{equation}
The value $2p^{\fr{k-1}{2}}\cdot \fr{\log \log p}{\sqrt{\log p}}$  is monotonically increasing in terms of $p$ and is larger than 5000, if $p\geq 5$.
Therefore, if $\tau(p)=\varepsilon t \ell, \ell \in L_{t}^{\varepsilon}$, we have $p\in S(\tau)$.

\section*{Acknowledgments}

The author would like to express his great gratitude to Masanobu Kaneko, who is my supervisor, Kyosuke Imanaka, Taiga Adachi and Keiichiro Nomoto for their comments and discussions with them.
This work was supported by JST SPRING, Japan Grant Number JPMJSP2136.

\section{Basic facts on arithemtic properties of $\tau$}

In this section, we recall the basic facts of $\tau$.
The following theorem was conjectured by Ramanujan, \eqref{RC1}, \eqref{RC2} was proved by Mordell \cite{Mordell}, \eqref{RC3} was proved by Deligne \cite{Deligne}.


\begin{theorem} \label{RCj}
\begin{enumerate}[(1)]
\item \label{RC1}
$\tau(nm)=\tau(n)\tau(m), \gcd(n, m)=1$,
\item \label{RC2}
$\tau(p^{m+1})=\tau(p)\tau(p^m)-p^{11}\tau(p^{m-1}), p:\text{prime}, m\geq 1$,
\item \label{RC3}
$|\tau(p)|<2p^{11/2}, p:\text{prime}.$
\end{enumerate}
\end{theorem}


Ramanujan showed that for $n\geq 1$ and prime $p\neq 23$,
\begin{equation} \label{RC}
\tau(n)\equiv 
\begin{cases}
n^2\sigma_1(n) \bmod 3, \\
n^3\sigma_1(n) \bmod 4, \\
n\sigma_1(n) \bmod 5, \\
n\sigma_3(n) \bmod 7,
\end{cases}
\tau(p)\equiv 0, 2, -1 \bmod 23.
\end{equation}



\begin{lemma} \label{lemmythmpf}
Let $p$ be prime.
\begin{enumerate}[(1)]
	\item \cite[Lemma 3.4]{BGPS} \label{lemfromBGPSlem3.4}
	We assume $\tau(p)\neq 0$.
	Let $r={\rm ord}_p(\tau(p))$, which means that $p$ can divide $\tau(p)$ exactly $r$ times.
	Then, for any $m\geq 1$, we have 
	$$
	{\rm ord}_p(\tau(p^m))=rm.
	$$

	\item \cite[Lemma 3.3]{BGPS} \label{lemfromBGPSlem3.3}
	If $\tau(p)=0$, we have
	$$
	\tau(p^m)=
	\begin{cases}
	0 & m\text{ : odd}, \\
	(-p^{11})^{m/2} & m\text{ : even}. 
	\end{cases}
	$$
\end{enumerate}
\end{lemma}

\begin{remark}
In Lemma \ref{lemmythmpf} \eqref{lemfromBGPSlem3.4}, $r$ is less than or equal to 5 from Deligne's bound $|\tau(p)|<2p^{11/2}$.
\end{remark}


\section{Tools}


To prove main results, many important tools are used.
In this section,  we introduce them and provide some propositions and lemmas.


\subsection{Lucas sequence} \label{Ls}


In this subsection, we introduce theory of Lucas sequence, which makes an important role to show that $d$ is small when
\begin{equation} \label{geneq}
\tau(p^{d-1})=a.
\end{equation}
(See Lemma \ref{methodEVT}).

For example, Balakrishnan, Ono, Craig, Tsai proved the following proposition by using the theory of Lucas sequences in \cite{BCKT}:
\begin{proposition}(\cite[Theorem 1.1]{BCKT}) \label{BCKT}
For any odd ordinary prime $\ell$, if $\tau(n)=\pm \ell^j, j\geq 1$, we must have $n=p^{d-1}$, where $p$ and $d$ are odd primes with $d \mid \ell(\ell^2-1)$.
\end{proposition}
Therefore, when $\tau(n)=\pm \ell^j$, we have $\omega(n)=1, \Omega(n)=d-1<\ell$.

The sequence $\{\tau(p^{i-1})\}_{i=1}^{\infty}$ forms a Lucas sequence in most cases.
A {\it Lucas pair} is a pair $(\alpha, \beta)$ of algebraic integers such that $\alpha+\beta, \alpha\beta\in \Z$ are nonzero coprime and $\alpha/\beta$ is not a root of unity.
Then, given a Lucas pair $(\alpha, \beta)$, we define
$$
u_n(\alpha, \beta):=\fr{\alpha^n-\beta^n}{\alpha-\beta}, n\geq 1
$$
and the sequence $\{u_n(\alpha, \beta)\}_{n=1}^{\infty}$ is called the {\it Lucas sequence} associated to the Lucas pair $(\alpha, \beta)$.

Let $(\alpha, \beta)$ be a Lucas pair and $\{u_n(\alpha, \beta)\}_{n=1}^{\infty}$ the Lucas sequence associated to it.
If  a prime $p \mid u_n$ does not divide $(\alpha-\beta)^2u_1\cdots u_{n-1}$, then $p$ is called {\it primitive prime divisor} of $u_n$.
For $n>2$, $u_n$ is {\it defective} if $u_n$ does not have any primitive prime divisor of $u_n$.

In the proof of Lemma \ref{methodEVT}, it is essential that Bilu, Hanrot, and Voutier \cite{Lucas} classified defective Lucas numbers completely.
They list up all defective Lucas numbers in Table 1 and 3 in \cite{Lucas}.
We will recall some classical facts on Lucas sequences to prove Lemma \ref{methodEVT}. 


\begin{lemma}[\cite{Lucas}, Proposition 2.1 (ii)] \label{Lucasgt}
For any Lucas sequence $\{u_n\}$, $u_d \mid u_n$ if $d \mid n$.
\end{lemma}


\begin{remark}
Proposition 2.1 (ii) in \cite{Lucas} is stated for Lehmer sequence.
On the other hand, if $d \mid n$, then, $(d, n)\not \equiv (0, 1) \bmod 2$.
For $n\geq 1$,
$$
u_n=
\begin{cases}
\widetilde{u}_n & n : \text{odd }, \\
(\alpha+\beta)\widetilde{u}_n=u_2\widetilde{u}_n & n : \text{even.}
\end{cases}
$$
Thus, the conclusion of this lemma is true.
\end{remark}


For any prime $\ell$, we put
$$
m_{\ell}(\alpha, \beta):=\min \{m\geq 1 : \ell \mid u_m\}.
$$
Note that we have $m_{\ell}(\alpha, \beta)>1$ by $u_1(\alpha, \beta)=1$.
For odd prime $\ell$ not divide $\alpha\beta$, $m_{\ell}(\alpha, \beta)<\infty$ by Corollary 2.2 in \cite{Lucas}.


\begin{proposition}[\cite{EVT}, Proposition 2.3] \label{rankell}
Let $(\alpha, \beta)$ be a Lucas pair and $\ell$ an odd prime such that $\ell$ does not divide $\alpha\beta(\alpha+\beta)$.
Then the following are true:
\begin{enumerate}[(1)]

	\item \label{discdivell}
	If $\ell \mid (\alpha-\beta)^2$, $m_{\ell}(\alpha, \beta)=\ell$.


	\item
	Otherwise, $m_{\ell}(\alpha, \beta) \mid (\ell-1)$ or $m_{\ell}(\alpha, \beta) \mid (\ell+1)$.

\end{enumerate}
\end{proposition}


We note that $\ell \mid (\alpha+\beta)$ is equivalent to $m_{\ell}(\alpha, \beta)=2$ by $u_2(\alpha, \beta)=\alpha+\beta$.


\begin{proposition}[\cite{Lucas} Corollary 2.2] \label{rkellap}
We assume that the notation is as Proposition \ref{rankell}.
	For any $m \geq 1$,
	\begin{equation}
	\ell \mid u_m(\alpha, \beta) \Longleftrightarrow m_{\ell}(\alpha, \beta) \mid m
\end{equation}
	holds.
\end{proposition}


\begin{remark}
In Corollary 2.2 of \cite{Lucas}, it is stated for Lehmer sequence.
However, under the assumption $\ell \nmid (\alpha+\beta)$, we have $\ell \mid u_m(\alpha, \beta)$ if and only if $\ell \mid \tilde{u}_m(\alpha, \beta)$.
Thus, $m_{\ell}$ in Corollary 2.2 of \cite{Lucas} is equal to $m_{\ell}$ in our notation.
\end{remark}


To show Lemma \ref{genEVTLemma2.1}, we will show the following lemma.


\begin{lemma} \label{DE}
Let $\ell$ and $p$ be odd primes, $\varepsilon=\pm 1$, $d\geq 3$ odd integer and $j, r\geq 0$ integer.

\begin{enumerate}[(1)]
\item \label{ei3decomp}
The only solutions of the equation
$$
p^2-\varepsilon p+1=3^r
$$
is $(\varepsilon, p, r)=(1, 2, 1)$.

\item \label{Lucassol3}
If
$$
2\ell^j=p^d+\varepsilon,
$$
we must have $p=3, \varepsilon=-1$.

\item \label{Lucassol7}
If
$$
8\ell^j=p^d+\varepsilon
$$
we must have $p=7, \varepsilon=1$.
\end{enumerate}

\end{lemma}


\begin{proof}
\eqref{ei3decomp}
We have $r>0$ immediately.
When $r=1$, by simple calculation, we have
$$
\varepsilon p (\varepsilon p-1)=2.
$$
Therefore, $p=2, \varepsilon=1$.

When $r\geq 2$, $\rho:=e^{2\pi i/3}$.
$$
(\varepsilon p+\rho)(\varepsilon p+\ol{\rho})=(1-\rho)^{2r}
$$
as equality of ideals in $\Z[\rho]$, where $\ol{(\cdot)}$ denote the complex conjugate.
There exists $u\in \Z[\rho]^{\times}(=\{\pm 1, \pm \rho, \pm \ol{\rho}\})$ such that $\varepsilon p+\rho=u(1-\rho)^r, \varepsilon p+\ol{\rho}=\ol{u}(1-\ol{\rho})^r$.
By $\varepsilon p \in \Z$, $\varepsilon p-\ol{\varepsilon p}=0$.
However, $\varepsilon p-\ol{\varepsilon p}=u(1-\rho)^r-\rho-\{\ol{u}(1-\ol{\rho})^r-\ol{\rho}\}\equiv -\rho(1-\rho) \bmod (1-\rho)^2$.
This is a contradiction noting $(1-\rho)$ is a prime ideal of $\Z[\rho]$.

\eqref{Lucassol3}, \eqref{Lucassol7}
The pair $(\alpha, \beta):=(p, -\varepsilon)$ is a Lucas pair.
Then, we can write $u_d(\alpha, \beta)=\fr{\alpha^d-\beta^d}{\alpha-\beta}=\fr{p^d+\varepsilon}{p+\varepsilon}$ since $d$ is odd, and 
\begin{equation} \label{Lucassol3eq}
(p+\varepsilon)u_d(\alpha, \beta)=2\ell^j, 8\ell^j.
\end{equation}
Since $\varepsilon=\pm 1$ is odd, we have $p\neq 2$.
Then, $p+\varepsilon$ is even, $u_d(\alpha, \beta)=\sum_{i=0}^{d-1}\alpha^{d-1-i}\beta^i$ is odd.
Thus, by \eqref{Lucassol3eq}, we can write
$p+\varepsilon=2\ell^s, 8\ell^s, u_d(\alpha, \beta)=\ell^{j-s}$, where $0\leq s \leq j$. 
If $s=j$, then $u_d(\alpha, \beta)=1$, this is impossible by $d\neq 1$.
Thus, we have $s<j$.
If $s>0$, then we have $\ell \mid (p+\varepsilon)^2$.
Thus, $u_d(\alpha, \beta)=\ell^{j-s}$ is defective.
Then, we must have $d=3$ by Table 1, 3 in \cite{Lucas} which tells us all defective Lucas numbers.
In fact, since $d$ is odd, we have $d=3, 5, 7, 13$.
However, since $b=(p+\varepsilon)^2$ is square, we have $d\neq 5, 7, 13$.
Since $\ell$ does not divide $p-\varepsilon$ (recall that $\ell $ divides $p+\varepsilon$), $m_{\ell}(\alpha, \beta)=3$.
By Proposition \ref{rankell} \eqref{discdivell}, $m_{\ell}(\alpha, \beta)=\ell$.
Thus, we have $\ell=3$.
Then, $u_3(p, -\varepsilon)=3^{j-s}$. and \eqref{ei3decomp} of this lemma implies $p=2$.
This is a contradiction.
Therefore, we have $s=0$.
In the case \eqref{Lucassol3}, $p+\varepsilon=2$.
We must have $p=3, \varepsilon=-1$.
In the case \eqref{Lucassol7}, we must have $p=7, \varepsilon=1$.
\end{proof}


\begin{lemma} \label{genEVTLemma2.1}
Let $k\geq 2$ be an integer and $A\in \Z$.
Suppose Lucas pair $(\alpha, \beta)$ are roots of the polynomial 
$$
X^2-AX+p^{2k-1},
$$
where $p$ is prime, $p \nmid A$ and $|A|\leq 2p^{(2k-1)/2}$.
Then for $n>2$, if $u_n(\alpha, \beta)\in \{\pm 4\ell^j, \pm 8\ell^j : j\geq 0, \ell : \text{odd prime}\}$, one of the following is true:
\begin{enumerate}[(1)]
\item \label{nondefective}
There exists a primitive prime divisor of $u_n$.
\item \label{defective}
We have $n=4$ and  
\begin{equation} \label{ell37Lucas}
\displaystyle 
\ell^{2j}=
\begin{cases}
\fr{3^{2k-1}-1}{2} & u_n=\pm 4\ell^j, \\
\fr{7^{2k-1}+1}{8}& u_n=\pm 8\ell^j.
\end{cases}
\end{equation}
\end{enumerate}
\end{lemma}


\begin{proof}
We assume that $u_n$ does not have any primitive prime divisor, in other words, $u_n$ is defective.
Immediately, it is showed that any element of $\{\pm 4\ell^j, \pm 8\ell^j : \ell \text{ odd prime}, j\geq 0\}$ never arises as defective number except for rows four of Table 2 in \cite{BCKT}. 
Then, we have $n=4$.
The information of the table implies that if $u_n=\pm 4\ell^j$, we have $m^2=2p^{2k-1}+\varepsilon 2$, where $m=2\ell^j, \varepsilon=\pm 1$.
Thus, we have
\begin{equation}
2\ell^{2j}=p^{2k-1}+\varepsilon.
\end{equation}
By Lemma \ref{DE} \eqref{Lucassol3}, $p=3, \varepsilon=-1$.
Therefore, \eqref{ell37Lucas} is true.

Similarly, if $u_n=\pm 8\ell^j$, we have
\begin{equation}
8\ell^{2j}=p^{2k-1}+\varepsilon.
\end{equation}
By Lemma \ref{DE} \eqref{Lucassol7}, $p=7, \varepsilon=1$.
Therefore, \eqref{ell37Lucas} is true.

\end{proof}


\begin{remark}
If $k=3, \ell=11, j=1$, then the first equation in \eqref{ell37Lucas} holds.
For $k<100000$, then the second equation in \eqref{ell37Lucas} holds.
Of course, if $(k, j)=(1, 0)$, for any $\ell$, two equations in \eqref{ell37Lucas} holds.
\end{remark}


\begin{proposition} \label{EVTLemma2.1}
Suppose $p$ is prime such that $p \nmid \tau(p)$, and $\alpha, \beta$ are roots of the integral polynomial
$$
X^2-\tau(p)X+p^{11}=(X-\alpha)(X-\beta).
$$
Then, we have $u_i(\alpha, \beta)=\tau(p^{i-1}), i\geq 1$ and the sequence $\{u_i(\alpha, \beta)\}$ satisfies following:
\begin{enumerate}[(1)]
\item (\cite[Lemma 2.1]{BCO}) \label{BCOlemma2.1}
For $n>2$, if $u_n(\alpha, \beta)\in \{\pm 1, \pm \ell : \ell \text{ odd prime}\}$, then $u_n(\alpha, \beta)$ is not defective.

\item (\cite[Lemma 2.1]{EVT}) 
For $n>2$, if $u_n(\alpha, \beta)\in \{\pm 2\ell^j : \ell \text{ odd prime}, j\geq 0\}$, then $u_n(\alpha, \beta)$ is not defective.

\item \label{mylem2.1BCO}
For $n>2$, if $u_n(\alpha, \beta)\in \{\pm 4\ell^j, \pm 8\ell^j : \ell \text{ odd prime}, j\geq 0\}$, then $u_n(\alpha, \beta)$ is not defective.
\end{enumerate}
\end{proposition}


\begin{proof}[proof of \eqref{mylem2.1BCO}]
If $2k-1=11$, the equation \eqref{ell37Lucas} has no solution.

\end{proof}


For $k\geq 1$, the equation $\tau(n)=2^k$ is still open, however, the following lemma is true:
\begin{lemma}[\cite{LL}] \label{pmpower2}
Let $k$ be an integer $\in \{1, \ldots, 6\}$.
The equation 
$$
\tau(n)=-2^k
$$
has no solution.
If $\tau(n)=2^k$, then we must have $n=p_1\cdots p_k$, where $p_1,\ldots, p_k$ are distinct primes such that $\tau(p_1)=\cdots =\tau(p_k)=2$.
\end{lemma}


This lemma is goal of this subsection.



\begin{lemma} \label{methodEVT}
Let $p, \ell$ be primes, $\ell \neq 2, d\geq 1$.
If
$$
\tau(p^{d-1})=\pm 2\ell, \pm 4\ell, \pm 8\ell,
$$
we have 
$$
d=
\begin{cases}
2, & \tau(p^{d-1})=\pm 2 \ell, \pm 8\ell, \\
2, 4, & \tau(p^{d-1})=\pm 4\ell.
\end{cases}
$$
Moreover, we have that $\{\tau(p^{i-1})\}_{i=1}^{\infty}$ is a Lucas sequence.
If $\tau(p^3)=\pm 4\ell$, then we have $\tau(p)=2$ and that $\ell$ is primitive prime divisor of $\tau(p^3)$.
\end{lemma}


\begin{proof}
We put $\tau(p^{d-1})=n$, where $n=\pm 2\ell, \pm 4\ell, \pm 8\ell$.
Clearly, $d>1$.

We will show that $d$ is even.
If $p=2$, by ${\rm ord}_2(\tau(2))={\rm ord}_2(-24)=3$ and Lemma \ref{lemmythmpf} \eqref{lemfromBGPSlem3.4}, we have $2^{3(d-1)} \mid \tau(2^{d-1})=n$.
We must have $d=2$.
Thus, we may assume that $p\neq 2$.
Since $n$ is even, \eqref{tauodd} implies that $d$ is even.

We will consider whether $p \nmid \tau(p)$ or not.
Lemma \ref{lemmythmpf} \eqref{lemfromBGPSlem3.3} implies $\tau(p)\neq 0$.
If $p \mid \tau(p)$, we have $p^{d-1} \mid \tau(p^{d-1})=n$ by Lemma \ref{lemmythmpf} \eqref{lemfromBGPSlem3.4}.
In this case, we must have $d=2$.
Therefore, we may assume that $p \nmid \tau(p)$, in other words, $\{1, \tau(p), \ldots \}$ is a Lucas sequence.
Let $u_i=\tau(p^{i-1})$ for $i\geq 1$.

We assume that $u_d$ does not have primitive prime divisor $\ell$.
If $d>2$, we have that $u_d$ does not have primitive prime divisor $2$ since we have $2 \mid u_2$ by \eqref{tauodd}.
This is a contradiction by Proposition \ref{EVTLemma2.1}.
Thus we have $d=2$.

If $u_d$ has primitive prime divisor $\ell$, we have $\ell \nmid u_1\cdots u_{d-1}$.
Since $d$ is even, we can put $m:=d/2$.
Then $\ell \nmid u_m$.
By Lemma \ref{Lucasgt}, $u_m \mid u_d$.
Therefore, $u_m=\pm 1, \pm 2, \pm 4, \pm 8$.
By Lemma \ref{pmpower2}, we must have $m=2, u_2=2$.
Thus $d=4$.
By the recurrence formula of $\tau(p^m)$, we have $\tau(p^3)=\tau(p)(\tau(p)^2-2p^{11})$.
By noting $u_2=\tau(p)$, we have $\tau(p^3)=4(2-p^{11})$.
Therefore, we have $u_4=\tau(p^3)\neq \pm 2\ell, \pm 8\ell$.
\end{proof}



\subsection{Thue equation}


For $a=\pm \ell$, we consider the equation \eqref{geneq}.
Thanks to Proposition \ref{BCKT}, there are only finitely many cases of possible values of $d$.
For each such $d$, applying the congruence \eqref{RC}, more explicitly, Lemma \ref{doddcongruences}, we can rule out many cases of the pair $(\ell, d)$.
When we cannot rule out $(\ell, d)$ by the congruences, we use the theory of Thue equation.

Let $F(X, Y)$ be a homogeneous polynomial with degree $\geq 3$ and $a\in \Z$.
An equation of the form
$$
F(X, Y)=a
$$
is called a {\it Thue equation}.

Comparing the Euler factor of the $L$-function of $\Delta$-function:
\begin{equation}
L(s, \Delta)=\sum_{n=1}^{\infty}\fr{\tau(n)}{n^s}=\prod_{p: \text{prime}}\fr{1}{1-\tau(p)p^{-s}+p^{11}p^{-2s}},
\end{equation}
we define 
\begin{equation}
\fr{1}{1-\sqrt{Y}T+XT^2}=\sum_{m=0}^{\infty}F_m(X, Y)T^m=1+\sqrt{Y}T+(Y-X)T^2+\cdots.
\end{equation}

Then, by simple calculation, for a nonnegative integer $m$, we have that $F_{2m}, F_{2m+1}/\sqrt{Y} \in \Z[X, Y]$ is a homogeneous polynomial with degree $m$ and $\tau(p^{2m})=F_{2m}(p^{11}, \tau(p)^2)$.
Some of the first of them are as follows:
\begin{gather}
F_0(X, Y)=1, F_1(X, Y)=\sqrt{Y}, F_2(X, Y)=Y-X, F_3(X, Y)=\sqrt{Y}(Y-2X), F_4(X, Y)=X^2-3XY+Y^2, \ldots
\end{gather}


\begin{proposition}[\cite{BCKT}, Lemma 5.1]
In particular, if $\tau(p^4)=a$, then $(x, y)=(p, 2\tau(p)^2-3p^{11})$ is an integer point on 
\begin{equation} \label{taup4eq}
y^2=5x^{22}+4a.
\end{equation}
\end{proposition}


Calculating by PARI/GP, for $d>5$, we can know whether the Thue equation
$$
F_{d-1}(X, Y)=\pm \ell
$$
has any solution $(X, Y)$ which can be written as $(X, Y)=(p^{11}, \tau(p)^2)$ or not.


\begin{lemma} \label{Thue}
Let $p$ be a prime.
If $(\ell, d)\in LD_1^{\varepsilon}$, then
$$
\tau(p^{d-1})\neq \varepsilon\ell,
$$
where
\begin{gather}
LD_1^{+}:=\{(277, 23), (421, 7), (631, 79), (827, 23), (827, 59), (967, 7), (967, 11), (967, 23)\}, \\
LD_1^{-}:=\{(367, 23), (443, 17), (643, 23), (643, 107), (827, 59), (829, 23), (829, 83), (919, 17)\}.
\end{gather}
\end{lemma}


\subsection{Dembner--Jain's method}


When $d=5$, the equation \eqref{geneq} for $a=281, 461, -599, -919$ cannot be ruled out by the congruences.
In this subsection, we will show that Dembner--Jain's method \cite{DJ} makes it possible ruling out $281, -919$.


\begin{lemma} \label{DJmethod}
For any prime $p$, we have $\tau(p^4)\neq 281, -919$.
\end{lemma}
\begin{proof}
We apply Dembner--Jain's method.
Let $\tau(p^4)=\varepsilon \ell$, where $(\varepsilon, \ell)=(1, 281), (-1, 919)$.
Let the real quadratic field $K:=\Q(\sqrt{5})$, $\mathcal{O}_K:=\Z[\omega]$ the ring of integers of $K$, where $\omega=(1+\sqrt{5})/2$.
Then, $\mathcal{O}_K^{\times}=\{\pm \omega^n \mid n\in \Z\}$.
In $K$, two prime ideal $281\Z, 919\Z$ is totally decomposition as follows: 
$$
(281)=(19+4\sqrt{5})(19-4\sqrt{5}), (919)=(42+13\sqrt{5})(42-13\sqrt{5}).
$$
When $\ell=281$, by \eqref{taup4eq},
$$
(19+4\sqrt{5})(19-4\sqrt{5})=\fr{y+p^{11}\sqrt{5}}{2}\cdot \fr{y-p^{11}\sqrt{5}}{2},
$$
where $y=2\tau(p)^2-3p^{11}$.
Let $\alpha=\fr{y+p^{11}\sqrt{5}}{2}$,
then there exists $u \in \mathcal{O}_K^{\times}$ such that
$$
\alpha=u(19+4\sqrt{5}) \text{ or } u(19-4\sqrt{5}).
$$
For simplicity, we assume the former case.
Then, taking its conjugate, we have 
$$
\ol{\alpha}=\fr{y-p^{11}\sqrt{5}}{2}=\ol{u}(19-4\sqrt{5}).
$$

We can write $u=\pm \omega^n, n\in \Z$.
Taking $\pm(\alpha-\ol{\alpha})$, we have
\begin{equation} \label{FT281}
(\pm p)^{11}=19\cdot \fr{\omega^n-\ol{\omega}^n}{\omega-\ol{\omega}}+4(\omega^n+\ol{\omega}^n)
\end{equation}
Here, we define the Fibonacci numbers and Lucas numbers
$$
F_n:=\fr{\omega^n-\ol{\omega}^n}{\omega-\ol{\omega}}, L_n:=\omega^n+\ol{\omega}^n
$$
(Note that "Lucas numbers" means specific numbers in this paper.
This is different to the notion of Lucas sequence.
In fact, the sequence of Lucas numbers $\{L_n\}_{n\geq 1}$ is not a Lucas sequence.)
The sequences $F_n, L_n \bmod 89$ has period 44 (we can see this from Debner--Jain's page https://github .com /sdembner /tauvalues), therefore, so is $x_n=19F_n+4L_n$.
There are no perfect 11-th power $\bmod 89$ in $x_0,\ldots, x_{43} \bmod 89$ by using computer.
Thus, the equation \eqref{FT281} has no solution.

When $\ell=919$, instead of \eqref{FT281},
\begin{equation} \label{FT919}
(\pm p)^{11}=42F_n+13L_n
\end{equation}

The sequence $x_n:=42F_n+13L_n$ has period 44 $\bmod 89$, 22 $\bmod 199$.
Then, $x_n$ is perfect 11-th power $\bmod 89$ if and only if $n\equiv 1, 12, 23, 34 \bmod 44$.
$x_n$ is perfect 11-th power $\bmod 199$ if and only if $n\equiv 11 \bmod 22$.
For any $n\in \Z$, these two condition can not be satisfied simultaneously.
Thus, the equation \eqref{FT919} has no solution.
\end{proof}


\section{Proof of main theorems}


In the last section, we will show the main theorems using the facts we have proved in the previous sections.


\subsection{proof of Theorem \ref{mainthmpmell}}

We let $\tau(n)=\varepsilon \ell, n\geq 1$.
By \eqref{lemfromLMcor1.3}, we may assume $\ell>252$.
Proposition \ref{BCKT} imply that $n=p^{d-1}$, where $d, p$ are odd primes such that $d \mid \ell(\ell^2-1)$.
Here, note that non-ordinary prime $\ell<1000$ is only $2, 3, 5, 7$(the next non-ordinary prime after 7 is 2411).

it is trivial that $\tau(p)\neq 0$ from Lemma \ref{lemmythmpf} \eqref{lemfromBGPSlem3.3}.
If $p \mid \tau(p)$, we have $p^{d-1} \mid \tau(p^{d-1})$ by Lemma \ref{lemmythmpf} \eqref{lemfromBGPSlem3.4}.
However, this is a contradiction by $d\geq 3$.
Therefore, we have that $p \nmid \tau(p)$.
Taking $\alpha, \beta$ as Proposition \ref{EVTLemma2.1}, the sequence $\{u_i(\alpha, \beta)=\tau(p^{i-1})\}_{i\geq 1}$ is a Lucas sequence.

Now, $u_d$ is not defective by Proposition \ref{EVTLemma2.1} \eqref{BCOlemma2.1} and $d\geq 3$.
In other words, $\ell$ is unique primitive prime divisor of $u_d$.
By definition, we have $\ell \nmid (\alpha-\beta)^2u_1\cdots u_{d-1}$.
In particular, $\ell \nmid (\alpha-\beta)^2$ and $\ell \nmid u_2=\alpha+\beta$.
We can check $\ell \nmid \alpha\beta(=p^{11})$ since we have $p\neq \ell$ by noting $p\nmid \tau(p)$.
Then, Proposition \ref{rankell} implies $d \mid (\ell^2-1)$.
By Lemma \ref{doddcongruences}, $\tau(p^{d-1})=\varepsilon \ell$ has no solution except for $(\ell, d)$ are as follows:
$$
(277, 23), (281, 5), (421, 7), (461, 5), (631, 79), (827, 23), (827, 59), (967, 7), (967, 11), (967, 23), \text{ if } \varepsilon=1,
$$
$$
(367, 23), (443, 17), (599, 5), (643, 23), (643, 107), (827 , 59), (829, 23), (829, 83), (919, 5), (919, 17),  \text{ if } \varepsilon=-1.
$$

When $d>5$, we can rule out $(\ell, d)$ by Lemma \ref{Thue}, $(281, 5), (919, 5)$ by Lemma \ref{DJmethod}.


\subsection{Proof of Theorem \ref{mainthm248}}


We may consider only the value of $\tau(p)$.
By Ramanujan's congruences, we can restrict the value of $\tau(p)$ to the element of $L_{t}^{\pm}, t=2, 4, 8$.
More precisely, this lemma is true:


\begin{lemma} \label{mylem}
We take an odd prime $\ell<1000$, and 
let $\varepsilon=\pm 1, t=2, 4, 8$.
\begin{enumerate}[(1)]
\item \label{ppm2ell}
For any prime $p$, if
$$
\tau(p)=\varepsilon t\ell,
$$
we must have $\ell\in L_{t}^{\varepsilon}$.

\item \label{p3pm48}
For any prime $p$,
$$
\tau(p^3)\neq \pm 4\ell.
$$
\end{enumerate}
\end{lemma}


\begin{proof}

	\eqref{ppm2ell}
For any $(\varepsilon, \ell) \not \in  L_{t}^{\varepsilon}$, $\tau(p)=\varepsilon t \ell$ never happen by the Ramanujan's congruence \eqref{RC}.
Here, note that when $t=2$, we may assume that $100<\ell$ by \eqref{EVTcor1.2}.

	\eqref{p3pm48}
If $\tau(p^3)=\varepsilon 4\ell, \varepsilon=\pm 1$, by Lemma \ref{methodEVT}, we have $\tau(p)=2$.
By the recurrence formula of $\tau$-function, we have $\tau(p^3)=\tau(p)(\tau(p)^2-2p^{11})$, thus $\varepsilon \ell=2-p^{11}$.
Now, $\ell<1000$ and $|2-p^{11}|\geq 2046$, therefore, this is a contradiction.
\end{proof}


\begin{proof}[proof of Theorem \ref{mainthm248}]
We assume that $\tau(n)=\pm 2^j\ell$ for an odd prime $\ell<1000$ and $j=1, 2, 3$.

We consider when $2\leq \omega(n)=:m$.
Then, we can write $n=p_1^{d_1-1}\cdots p_m^{d_m-1}$, where $p_1,\ldots p_m$ are distinct primes and $d_1, \ldots, d_m\geq 2$.
We put $q_i:=p_i^{d_i-1}, i=1,\ldots, m$.
From the mulitiplicity of $\tau$, we have
$$
\prod_{i=1}^m|\tau(q_i)|=2^j\ell.
$$
Then, there exists $i$ such that $\ell \nmid |\tau(q_i)|$ because $m\geq 2$.
Then $|\tau(q_i)|=1, 2, \ldots, 2^j$.
By \eqref{tauneqpm1prime691}, $|\tau(q_i)|\neq 1$ and by Lemma \ref{pmpower2}, $d_i=2$ and $\tau(q_i)=2$.
Therefore, $n/q_i(=n/p_i)$ is coprime to $q_i(=p_i)$.
The case \eqref{fcmt} of Theorem \ref{mainthm248} is true.

If $\omega(n)=1$, so that we can write $n=p^{d-1}$, where $p$ is prime, $d\geq 2$ is an integer.
By Lemma \ref{methodEVT} and Lemma \ref{mylem} \eqref{p3pm48}, we have $d=2$.
By Lemma \ref{mylem} \eqref{ppm2ell}, \eqref{p3pm48}, we have $\ell \in L_{2^j}^{\varepsilon}$.
Thus, the case \eqref{scmt} of Theorem \ref{mainthm248} is true.
\end{proof}


\appendix

\section{Ramanujan's congruences for $\tau$-values at prime powers}
\begin{lemma} \label{doddcongruences}
For any odd $d$, we have
$$
\tau(p^{d-1})\underset{(12)}{\equiv}0, 1, d
$$
$$
\tau(p^{d-1})\underset{(5)}{\equiv}
\begin{cases}
0, 1, d & d \underset{(4)}{\equiv} 1, \\
0, 1, 2, 3, d & d \underset{(4)}{\equiv} 3.
\end{cases}
$$
$$
\tau(p^{d-1})\underset{(7)}{\equiv}
\begin{cases}
0, 1, d & d \underset{(6)}{\equiv} 1, \\
0, 1, 2, 4, d, 2d, 4d & d \underset{(6)}{\equiv} 3, 5.
\end{cases}
$$
$$
\tau(p^{d-1})\underset{(23)}{\equiv}
\begin{cases}
1, d & d \underset{(6)}{\equiv} 1, \\
0, 1, d& d \underset{(6)}{\equiv} 3, \\
1, -1, d & d \underset{(6)}{\equiv} 5
\end{cases}
$$
for any prime $p$.
\end{lemma}

\begin{lemma} \label{devencongruences}
If $d$ is even $\geq 2$, we have
$$
\tau(p^{d-1})\underset{(12)}{\equiv}0, d,
$$
$$
\tau(p^{d-1})\underset{(5)}{\equiv}
\begin{cases}
0, d & d \underset{(4)}{\equiv} 0, \\
0, 1, 2, d & d \underset{(4)}{\equiv} 2,
\end{cases}
$$
$$
\tau(p^{d-1})\underset{(7)}{\equiv}
\begin{cases}
0, d, 2d, 4d & d \underset{(6)}{\equiv} 0, 2, \\
0, d & d \underset{(6)}{\equiv} 4,
\end{cases}
$$
for any prime $p$ and
$$
\tau(p^{d-1})\underset{(23)}{\equiv}
\begin{cases}
0, d & d \underset{(6)}{\equiv} 0, \\
0, -1, d & d \underset{(6)}{\equiv} 2, \\
0, 1, d& d \underset{(6)}{\equiv} 4,
\end{cases}
$$
for any prime $p\neq 23$.
When $p=23$, $\tau(p^{d-1})\underset{(23)}{\equiv} 1$ for any integer $d\geq 1$.
\end{lemma}




\bibliographystyle{plain}
\bibliography{PDreference}

\begin{thebibliography}{10}

\bibitem{BCO}
Jennifer~S. Balakrishnan, William Craig, and Ken Ono.
\newblock Variations of {L}ehmer's conjecture for {R}amanujan's tau-function.
\newblock {\em J. Number Theory}, 237:3--14, 2022.

\bibitem{BCKT}
Jennifer~S. Balakrishnan, William Craig, Ken Ono, and Wei-Lun Tsai.
\newblock Variants of {L}ehmer's speculation for newforms.
\newblock {\em Adv. Math.}, 428:Paper No. 109141, 31, 2023.

\bibitem{EVT}
Jennifer~S. Balakrishnan, Ken Ono, and Wei-Lun Tsai.
\newblock Even values of {R}amanujan's tau-function.
\newblock {\em Matematica}, 1(2):395--403, 2022.

\bibitem{STS2}
Tom Barnet-Lamb, David Geraghty, Michael Harris, and Richard Taylor.
\newblock A family of {C}alabi-{Y}au varieties and potential automorphy {II}.
\newblock {\em Publ. Res. Inst. Math. Sci.}, 47(1):29--98, 2011.

\bibitem{BGPS}
Michael~A. Bennett, Adela Gherga, Vandita Patel, and Samir Siksek.
\newblock Odd values of the {R}amanujan tau function.
\newblock {\em Math. Ann.}, 382(1-2):203--238, 2022.

\bibitem{Lucas}
Yu. Bilu, G.~Hanrot, and P.~M. Voutier.
\newblock Existence of primitive divisors of {L}ucas and {L}ehmer numbers.
\newblock {\em J. Reine Angew. Math.}, 539:75--122, 2001.

\bibitem{Deligne}
Pierre Deligne.
\newblock La conjecture de {W}eil. {I}.
\newblock {\em Inst. Hautes \'{E}tudes Sci. Publ. Math.}, (43):273--307, 1974.

\bibitem{DJ}
Spencer Dembner and Vanshika Jain.
\newblock Hyperelliptic curves and newform coefficients.
\newblock {\em J. Number Theory}, 225:214--239, 2021.

\bibitem{GTW}
Ayla Gafni, Jesse Thorner, and Peng-Jie Wong.
\newblock Almost all primes satisfy the {A}tkin-{S}erre conjecture and are not
  extremal.
\newblock {\em Res. Number Theory}, 7(2):Paper No. 31, 5, 2021.

\bibitem{HM}
Mitsuki Hanada and Rachana Madhukara.
\newblock Fourier coefficients of level 1 {H}ecke eigenforms.
\newblock {\em Acta Arith.}, 200(4):371--388, 2021.

\bibitem{STS1}
Michael Harris, Nick Shepherd-Barron, and Richard Taylor.
\newblock A family of {C}alabi-{Y}au varieties and potential automorphy.
\newblock {\em Ann. of Math. (2)}, 171(2):779--813, 2010.

\bibitem{LL}
Kaya Lakein and Anne Larsen.
\newblock Some remarks on small values of {$\tau(n)$}.
\newblock {\em Arch. Math. (Basel)}, 117(6):635--645, 2021.

\bibitem{LT}
Serge Lang and Hale Trotter.
\newblock {\em Frobenius distributions in {${\rm GL}\sb{2}$}-extensions},
  volume Vol. 504 of {\em Lecture Notes in Mathematics}.
\newblock Springer-Verlag, Berlin-New York, 1976.
\newblock Distribution of Frobenius automorphisms in ${\rm
  GL}\sb{2}$-extensions of the rational numbers.

\bibitem{L}
D.~H. Lehmer.
\newblock The vanishing of {R}amanujan's function {$\tau(n)$}.
\newblock {\em Duke Math. J.}, 14:429--433, 1947.

\bibitem{LM}
Wenwen Lin and Wenjun Ma.
\newblock On values of {R}amanujan's tau function involving two prime factors.
\newblock {\em Ramanujan J.}, 63(1):131--155, 2024.

\bibitem{LR}
Nik Lygeros and Olivier Rozier.
\newblock Odd prime values of the {R}amanujan tau function.
\newblock {\em Ramanujan J.}, 32(2):269--280, 2013.

\bibitem{Mordell}
L.J. Mordell.
\newblock On mr. ramanujan's empirical expansions of modular functions.
\newblock {\em Proc. Camb. Philos. Soc.}, (19):117--124, 1917.

\bibitem{MMS}
M.~Ram Murty, V.~Kumar Murty, and T.~N. Shorey.
\newblock Odd values of the {R}amanujan {$\tau$}-function.
\newblock {\em Bull. Soc. Math. France}, 115(3):391--395, 1987.

\bibitem{NT1}
James Newton and Jack~A. Thorne.
\newblock Symmetric power functoriality for holomorphic modular forms, {I},
  {II}.
\newblock {\em Publ. Math. Inst. Hautes \'{E}tudes Sci.}, 134:1--152, 2021.

\bibitem{Serre}
Jean-Pierre Serre.
\newblock Divisibilit\'{e} de certaines fonctions arithm\'{e}tiques.
\newblock {\em Enseign. Math. (2)}, 22(3-4):227--260, 1976.

\bibitem{Thorner}
Jesse Thorner.
\newblock Effective forms of the {S}ato-{T}ate conjecture.
\newblock {\em Res. Math. Sci.}, 8(1):Paper No. 4, 21, 2021.

\end{thebibliography}


\end{document}